\documentclass[10pt]{article}
\usepackage{amssymb,amsmath,amsfonts}
\usepackage{dsfont}
\usepackage{amsthm}
\usepackage{graphicx,color,subfig}
\usepackage{xspace}

\usepackage[latin1]{inputenc}
\usepackage[T1]{fontenc}

\usepackage{mathrsfs}
\usepackage{amscd}

\DeclareMathOperator{\Tr}{Tr}

\newcommand{\Con}{\ensuremath{\mathscr C}}

\newcommand{\SSS}{\ensuremath{\mathscr S}}

\DeclareMathOperator{\supp}{supp}

\newcommand{\est}[1]{\langle #1 \rangle}

\newcommand{\mb}[1]{\ensuremath{\mathbb{#1}}}
\newcommand{\N}{{\mb{N}}}

\newcommand{\R}{{\mb{R}}}
\newcommand{\C}{{\mb{C}}}

\newcommand{\eps}{\varepsilon}

\renewcommand{\d}{\ensuremath{\partial}}

\newenvironment{prooff}[1][Proof]{\textbf{#1.} }{\hfill \rule{0.5em}{0.5em}}

\DeclareMathOperator{\Op}{Op}

\newcommand{\nhd}{neighborhood\xspace}

\usepackage[active]{srcltx}

\newcommand{\tno}[1]{ |\hskip -0,5pt |\hskip -0,5pt | #1 |\hskip -0,5pt |\hskip -0,5pt |}

\usepackage[english]{babel}

\newtheorem{lemma}{Lemma}[section]
\newtheorem{theorem}{Theorem}

\theoremstyle{definition}

\newtheorem{remark}{Remark}

\makeatletter
\def\keywords{
    \vspace{1ex}
    \noindent
    \if@twocolumn
      \small{\bf  Keywords}\/---$\!$    \else
      \begin{center}\small\ {\bf Keywords}\end{center}\quotation\small
    \fi}
\def\endkeywords{\vspace{0.6em}\par\if@twocolumn\else\endquotation\fi
    \normalsize\rm}
\makeatother

\begin{document}

\title{Counting function for interior transmission eigenvalues}

\author{ Luc Robbiano\thanks{Laboratoire de Math\'ematiques de Versailles, Universit\'e de Versailles St Quentin,
CNRS UMR 8100, 45, Avenue des \'Etats-Unis, 78035 Versailles, France. e-mail : luc.robbiano@uvsq.fr
}}

\maketitle
\begin{abstract}
In this paper we give results on the counting function associated with the interior transmission eigenvalues. For a complex refraction index we estimate of the counting function by $Ct^{n}$. In the case where the refraction index is positive we give an equivalent of the counting function. 
\end{abstract}


\tableofcontents

\begin{section}{Introduction}
In this paper we give an estimate to the counting function associated with the interior transmission eigenvalues. We recall the problem. Let $\Omega$ be a smooth bounded domain in $\R^n$. Let $n(x)$ be  a smooth function defined in $\overline \Omega$, called the refraction index. We say that $k\not=0$ is a interior transmission eigenvalue if there exists $(w,v)\not=(0,0)$ such that
\begin{equation}\label{eq : ite}
 \left\{
	\begin{array}{l}
	 	\Delta w +k^2n(x) w=0 \text{ in } \Omega,\\ 
		\Delta v +k^2 v=0 \text{ in } \Omega,\\
		w=v \text{ on } \d \Omega,\\
		\d_\nu w=\d_\nu v  \text{ on } \d \Omega,
	\end{array}\right.
\end{equation}
where $\d_\nu$ is the exterior normal derivative to $\d\Omega$.
We consider here the function $n(x)$ complex valued. In  physical models, we have $n(x)=n_1(x)+in_2(x)/k$ where $n_j$ are real valued. Taking $u=w-v$ and $\tilde v=k^2v$, we obtain the following equivalent system if $k\not=0$,
\begin{equation}\label{eq : ite Syl} 
\left\lbrace 
 	\begin{array}{ll}
 		&\big( \Delta +k^2(1+m)\big) u+mv=0\text{ in } \Omega,\\
		&(\Delta +k^2)v=0 \text{ in } \Omega,\\
		&u=\partial_\nu u=0 \text{ on } \partial \Omega,
	\end{array}\right. 
\end{equation}
where, for simplicity, we have replaced $\tilde v$ by $v$ and $n$ by $1+m$.

When $k\in\R$, this problem is related with scattering problem. We can find a precise result in  Colton and 
Kress~\cite[Theorem~8.9]{CK92} first proved by Colton, Kirsch and Päivärinta~\cite{CKP89} and in a survey  by Cakoni and Haddar~\cite{CH12}.

As the problem si not self-adjoint even for $n(x)$ real valued, usual tools used in self-adjoint cases cannot be applied, in particular, even for operator with compact resolvent, the existence of $k$ is not always true.

A lot of results was obtained this last years using several methods. 
When $n(x) $ is real, Päivärinta and Sylvester~\cite{PS08}  proved that there exist interior transmission eigenvalues; 
 Cakoni, Gintides, and Haddar~\cite{CGH10} proved that the set of $k^2_j$ is infinite and discrete. For $n(x)$ complex valued  Sylvester~\cite{Sy12} proved that this set is discrete finite or infinite.
  In~\cite{Ro2013} we proved that there exist a infinite number of complex eigenvalues and the associated generalized eigenspaces span a dense space in $L^{2}(\Omega)\oplus L^{2}(\Omega)$.

Other related problems are studied in literature, problem with cavities studied by Cakoni, Çayö\-ren, and Colton~\cite{CCC08},  Cakoni, Colton, and Haddar~\cite{CCH10}, problem for operators of order $m>2$
 by Hitrik, Krupchyk, Ola,  and Päivärinta~\cite{HKOP10, HKOP11}.
 
Lakshtanov and Vainberg~\cite{LV12, LV11-arxiv, LV12-arxiv} studied the counting function for problems with different boundary conditions. For Problem~\eqref{eq : ite},  in~\cite{LV12bis-arxiv}   they obtain a lower estimate of the counting function for real interior transmission eigenvalues.

For counting function in~\cite{Ro2013} we gave some non optimal estimate. This estimate was improved by 
Dimassi and Petkov~\cite{DP13-arx}.  Their estimate have the same size than the one found below in 
Theorem~\ref{th. count. funct.} except for  a factor $3\sqrt{3}$.  For constant $m$ there is a recent result given by Pham and Stefanov~\cite{PS2013} where they give an equivalent of counting function in this case.

The main results of this paper are Theorem~\ref{th. trace formula} and Theorem~\ref{th. count. funct.}. In 
Theorem~\ref{th. trace formula} we prove than the counting function satisfies an estimate in $Ct^{n}$ and we prove that if we denote by $\lambda_{j}$ the eigenvalues of the problem 
\begin{equation*}
\sum_{j\in\N} \frac{1}{\lambda_j^p-z^p}= (2\pi)^{-n} |z|^{-p+n/2}\!\!\!\int_\Omega\int\!\!\!\left(((1+m(x)^{-p}|\xi|^{2p}-\mu^p)^{-1} +(|\xi|^{2p}-\mu^p)^{-1}
\right) d\xi dx+o(|z|^{-p+n/2}),
\end{equation*}
when $|z|$ goes to $\infty$ and $z $ in a line outside a domain related with the range of $n(x)$.

In the case where $n(x)>0$ is real valued or  if $n_{1}(x)>0$ in the case $n(x)=n_{1}(x)+n_{2}(x)/k$ this estimate allows, applying a tauberian theorem to give an equivalent of the counting function. We find that $N(t)\sim \alpha t^{n}$ where the precise value of $\alpha$ is given in Theorem~\ref{th. count. funct.}. These results are also proven by Faierman~\cite{Fa2013} in a preprint. The methods used are very close that the one used here but he assumes that $n(x)\not=1$ every where. This condition excludes the case of cavity. Here we assume only that $n(x)\not=1$  in a \nhd of the boundary.

\end{section}

\begin{section}{Notations and background}

Let $\Omega $ be a $\Con^\infty$ bounded domain in $\R^n$. Let $n(x)\in \Con^\infty(\overline\Omega)$ be complex valued. We set  $m(x)=n(x)-1$. We consider also the case where $n(x)=n_1(x)+in_2(x)/k$ where $n_j(x)$ are real valued and $k$ the spectral parameter. This case is different of the previous one but can be treated similarly. We assume that for all $x\in\overline \Omega$, $n(x)\not=0$, or $n_1(x)\not=0$ or equivalently $m(x)\not= -1$. We assume that  there exists a \nhd $W$ of $\d \Omega$ such that for $x\in\overline W$,  $n(x)\not=1$ or  $n_1(x)\not=1$ or equivalently $m(x)\not= 0$. Actually if $n(x)\not=1$ for all $x\in\d\Omega$, such a \nhd $W$ exists.

We denote by 
 $C_e$   the cone in $\C$ defined by 
\begin{equation}\label{def C}
 C_e=\{ z\in\C,\ \exists x\in\overline \Omega,\ \exists \lambda\ge 0,\text{ such that }z=\lambda(1+\overline m(x))\}.
\end{equation}

In the case where $n(x)=n_1(x)+in_2(x)/k$, $C_e=[0,\infty)$ if $n_1(x)>0$ for all $x\in\overline \Omega$, and $C_e=[-\infty,0]$ if $n_1(x)<0$ for all $x\in\overline \Omega$.

Here we give some notations useful for the statement of the results. We use the notations and the results proven
in~\cite{Ro2013}, except some change of sign.

Let $z\in\C$, we denote by $B_z(u,v)=(f,g)$ the mapping defined from $H^2_0(\Omega)\oplus \{v\in L^2(\Omega), \ \Delta v\in L^2(\Omega)\}$ to $L^2(\Omega)\oplus L^2(\Omega)$  by

\begin{equation}\label{eq : resolvent 1}
\left\lbrace 
 	\begin{array}{ll}
 		&\big( \frac{-1}{1+m}\Delta -z\big) u-\frac{m}{1+m}v=f\text{ in } \Omega\\
		&(-\Delta -z)v=g \text{ in } \Omega\\
	\end{array}\right. 
\end{equation}
 
In the case where $n=n_1+in_2/k$ we must change the definition of $B_z$. We define $m_1(x)=n_1(x)-1$ and $m_2(x)=n_2(x)$. The mapping $\hat B_k(u,v)=(f,g)$ is given by
\begin{equation}\label{eq : resolvent 1 bis}
\left\lbrace 
 	\begin{array}{ll}
 		&\big( \frac{-1}{1+m_1}\Delta -k^2-ik\frac{m_2}{1+m_1}\big) u+\left( \frac{-m_1}{1+m_1}-\frac{im_2}{k(1+m_1)}\right) 
		v=f\text{ in } \Omega\\
		&(-\Delta -k^2)v=g \text{ in } \Omega\\
	\end{array}\right. 
\end{equation}

Remark  that the principal symbol of $\hat B_k$ is 
the same than the one of $B_z$ if we set $z=k^{2}$.

Under an assumption on $z$, $B_z$ is invertible for some $z$.

\begin{theorem}\label{th : Existence}
 Assume $C_e\not=\C$, then there exists $z\in \C$ such that $B_{z}$ is a  bijective  map from $H^2_0(\Omega)\oplus \{v\in L^2(\Omega), \ \Delta v\in L^2(\Omega)\}$ to $ L^2(\Omega)\oplus  L^2(\Omega)$.
\end{theorem}

In the case $n(x)=n_1(x)+n_2(x)/k$, here $C_e\not=\C$ and we have the same result.
\begin{theorem}\label{th : Existence bis}
There exists $k\in \C$ such that $\hat B_k$ is  bijective from $H^2_0(\Omega)\oplus \{v\in L^2(\Omega), \ \Delta v\in L^2(\Omega)\}$ to $ L^2(\Omega)\oplus  L^2(\Omega)$.
\end{theorem}

If for $z\in \C$ the solution $B_z(u,v)=(f,g)$ exists we  denote by $R_z(f,g)=(u,v)$. In case where $n(x)=n_1(x)+in_2(x)/k$ if for $k\in \C$ the solution of $\hat B_k(u,v)=(f,g)$ exists we denote by $\hat R_k(f,g)=(u,v)$.

\begin{theorem}\label{Resolv comp}
 Assume $C_e\not=\C$, there exists $z\in\C$ such that the resolvent $R_{z}$ from $\overline H^2(\Omega)\oplus L^2(\Omega)$ to itself is compact. 

In particular,  applying  the Riesz theory, the spectrum is finite or is a discrete countable set. If $\lambda\not=0$ is in the spectrum, $\lambda $ is an eigenvalue associated with a finite  dimensional  generalized eigenspace. 
\end{theorem}
 
\begin{theorem}\label{Resolv comp bis}
There exists $k\in \C$ such that the resolvent $\hat R_{k}$ from $\overline H^2(\Omega)\oplus L^2(\Omega)$ to itself is compact. 

In particular, we can apply the Riesz theory, the spectrum is finite or a discrete countable set. If $\lambda\not=0$ is in the spectrum, $\lambda $ is an eigenvalue associated with a finite  dimensional  generalized eigenspace. 
\end{theorem}

\begin{remark}\label{rem vp pres du reel}
 Actually if $z_0 \not\in C_e\cup [0,\infty)$ for all $\lambda>0 $ large enough we can take $z=\lambda z_0$ in  Theorems~\ref{th : Existence} and \ref{Resolv comp}.
 
 If $k_0^2\not\in C_e \cup [0,\infty)$  for all $\lambda>0 $ large enough we can take $k=\lambda k_0 $ in the Theorems~\ref{th : Existence bis} and \ref{Resolv comp bis}.  Here we estimate the resolvent in the exterior of a conic \nhd of $C_e \cup [0,\infty)$. In particular if $n_{1}(x)>0$, the eigenvalues $k^2$ are in all small conic \nhd of $(0,+\infty)$, except for  a finite number of eigenvalues.
  \end{remark}
 
 In general for a non self-adjoint problem,  we cannot claim that the spectrum is non empty. In the following theorem, with a stronger assumption on $C_e$, we can prove that the spectrum is non empty.

We say that $C_e$ is contained in a sector with angle less than $\theta$ if there exist $\theta_1<\theta_2$, such that $C_e\subset \{ z\in \C, \ z=0 \text{ or } \frac{z}{|z|}=e^{i\varphi},\text{ where }\theta_1\le \varphi\le\theta_2\}$, and $\theta_2-\theta_1\le\theta $.

\begin{theorem}\label{Spect infini}
 Assume that $C_e$ is contained in a sector with angle less than $\theta$ with $\theta<{2\pi}/{p}$  where $4p>n$ and $\theta<\pi/2$. Then there exists $z$ such that the spectrum of $R_{z}$ is infinite and the space spanned by the generalized eigenspaces is dense in $H^2_0(\Omega)\oplus\{ v\in L^2(\Omega),\ \Delta v\in L^2(\Omega)\}$. 
 \end{theorem}
 
 \begin{theorem}\label{Spect infini bis}
There exists $k$ such that the spectrum of $\hat R_{k}$ is infinite and the space spanned by the generalized eigenspaces is dense in $H^2_0(\Omega)\oplus\{ v\in L^2(\Omega),\ \Delta v\in L^2(\Omega)\}$. 
 \end{theorem}

\begin{remark}
 These results are based on the theory given in Agmon~\cite{Ag65} and using the spectral results on Hilbert-Schmidt operators. In this theory we deduce that the spectrum is infinite from the proof that the generalized eigenspaces form a dense  subspace in the closure of the range of $R_{z}$ [resp. $\hat R_k$]. In~\cite{Ro2013} we proved that $R_{z}^p$ [resp.$\hat R_k^p$ ] is a Hilbert-Schmidt operator if $4p>n$. We can deduce the spectral decomposition of $R_{z}$   [resp. $\hat R_k$] from that of $R_{z}^p$ [resp.$\hat R_k^p$ ].
\end{remark}

Let $z_j$  be the elements of the spectrum of $R_{z}$  [resp. $\hat R_k$]   and $E_j$ the generalized associated eigenspace. We denote by $N(t)=\sum_{|z_j|^{-1}\le t^2}\dim E_j$.

If $z_{j}$ is an eigenvalue of $R_{z}$, $\lambda_{j}=-z+1/z_{j}$ is an eigenvalue of $B_{0}$   
and we have $N(t)\sim \sharp\{j,\ |\lambda_{j} |\le t^{2}  \}$, where $\lambda_{j}$ are counted with multiplicity.

\end{section}

\begin{section}{Results}

We denote by $\omega_j$ for $j=1,\cdots, p$, the roots of $z^p=1$.

\begin{theorem}\label{th. trace formula} 

We assume as in Theorem~\ref{Spect infini} that $\theta<2\pi/p$ and $\theta<\pi/2$ where $p$ satisfies $2p>n$ and $4p>4+n$.
Then, there exists $C>0$ such that $N(t)\le Ct^{n}$. 

Moreover
 let $\mu\in\C$ such that $|\mu|=1$ and we assume that $\omega_{j}\mu\not\in C_{e}\cup (0,+\infty)$  for $j=1,\cdots, p$. We denote by $a(x)=(1+m(x))^{-1}$ or $a(x)=(1+m_{1}(x))^{-1}$ if $n(x)=n_{1} (x)+k^{-1}n_{2}(x)$. We fix $z_0$ such that the resolvent $R_{z_0}$ exist and let $\mu_j$ such that $1/\mu_j$ are the eigenvalues of $R_{z_0}$ counted with multiplicity. Then we have
\begin{equation}\label{form entre serie et integrale}
\sum_{j\in\N} \frac{1}{\mu_j^p-z^p}= (2\pi)^{-n} |z|^{-p+n/2}\int_\Omega\int\left((a^p|\xi|^{2p}-\mu^p)^{-1} +(|\xi|^{2p}-\mu^p)^{-1}
\right) d\xi dx+o(|z|^{-p+n/2}),
\end{equation}
when $z=r\mu$ and $r$ goes to $\infty$.
\end{theorem}

\begin{remark}
The first part of the theorem improve~\cite[Theorem~7]{Ro2013} where we  found  the estimate $N(t)\le Ct^{n+4}$.
\end{remark}

\begin{theorem}\label{th. count. funct.}
We assume 
$n_{1}(x)=1+m_{1}(x)>0$ for all $x\in\overline\Omega$ when $n(x)=n_{1}(x)+n_{2}/k$. Then
$$N(t)= \alpha t^n+o(t^n) \text{ where }\alpha =(2\pi)^{-n}\text{Vol}(B_1) \int_\Omega((1+m_1(x))^{n/2}(x) +1)dx.$$
\end{theorem}

\begin{remark} 
By a more precise study in a \nhd of the boundary we can obtain a result with a smaller remainder 
in~\eqref{form entre serie et integrale} but this estimate does not allow to prove a result better on the counting function. Malliavin~\cite{Ma1962} was proved a tauberian  theorem with a sharp remainder but this requires an estimate on $\sum_{j\in\N} \frac{1}{\lambda_j^p-z^p}$ in a complex domain except in a parabola \nhd of $(0,\infty)$. Here the estimate are proved in a complex domain except in a conic \nhd of $(0,\infty)$.  It is maybe possible to improve this result following  Hitrik, Krupchyk, Ola and Päivärinta~\cite{HKOP-arx} where they prove  that the eigenvalues are in a parabolic  \nhd of $(0,\infty)$.
\end{remark}

\end{section}

\begin{section}{Proof of Theorem \ref{th. trace formula}}

The proof is based on Lemmas~\ref{Lem. Kernel trace} and \ref{Lem. trace symbol} below. We introduce some notations.

We set $S=R_{z_{0}}$ 
and $T=S^p$ where $p$ satisfies the assumption of Theorem~\ref{th. trace formula}. 

We set $T_\lambda= T(I-\lambda T)^{-1}$.  

As $T_\lambda $ is a matrix of operators we denote
$$
T_\lambda=
\begin{pmatrix}
T_{11} & T_{12}\\ 
T_{21}& T_{22}
\end{pmatrix},
$$

We denote by 
\begin{equation}\label{Def. noyau Tzp}
\begin{pmatrix}
K_{11} &K_{12}  \\ K_{21}  &K_{22}
\end{pmatrix},
\end{equation}
the kernel of $T_{z^p}$.

\begin{lemma}\label{Lem. Kernel trace}

Under the assumption of Theorem~\ref{th. trace formula},
there exists $C>0$ such that $N(t)\le Ct^{n}$. 

Moreover let $V$ a conical \nhd of $C_{e} \cup [0,\infty)$, there exists $R>0$ such all $z\in\C$, satisfying $|z|\ge R$ and $\omega_{j}z\notin\overline V$, we have
\begin{equation}
\int_\Omega K_{11}(x,x)dx+\int_\Omega K_{22}(x,x)dx=\sum_{j\in\N} \frac{1}{\lambda_j^p-z^p}. 
\end{equation}
\end{lemma}

\begin{lemma}    \label{Lem. trace symbol} 
With the notation of Theorem~\ref{th. trace formula}
we have
$ |z|^{p-n/2}
\left( \int_\Omega K_{11}(x,x)dx+\int_\Omega K_{22}(x,x)dx\right)$   
goes to $(2\pi)^{-n}\int_\Omega\int\left((a^p|\xi|^{2p}-\mu^p)^{-1} +(|\xi|^{2p}-\mu^p)^{-1}
\right) d\xi dx$ when $|z|$ goes to $\infty$.
\end{lemma}

Clearly Lemmas~\ref{Lem. Kernel trace} and \ref{Lem. trace symbol} imply Theorem~\ref{th. trace formula}.

\begin{subsection}{Proof of Lemma \ref{Lem. Kernel trace}} 

We recall that  $\omega_j$ for $j=1,\cdots, p$, are the roots of $z^p=1$, we have
\begin{equation}
( 1-z^pS^p)=\prod_{j=1}^p(1-\omega_jzS).
\end{equation}
The operator $  ( 1-z^pS^p) $ is invertible if and only if $  (1-\omega_jzS) $ is invertible  for all $j$. Thus we have 
$$
T_{z^p}=S^p\prod_{j=1}^p(1-\omega_jzS)^{-1}=\prod_{j=1}^p S_{\omega_jz}.
$$

If we denote by
$$
S_z=
\begin{pmatrix}
S_{11} &S_{12} \\ S_{21}& S_{22}
\end{pmatrix}.
$$

We  recall that $S=R_{z_{0}} $ and 
\begin{equation} \label{Lien S et R}
S_z=(R_{z_{0}})_z=R_{z_{0}}(I-z R_{z_{0}})^{-1}=(R_{z_{0}}^{-1}-z)^{-1}=(B_{z_{0}}-z)^{-1}=(B_0-z_{0}-z)^{-1}=R_{z_{0}+z},
\end{equation}
if the resolvent $R_{z_{0}+z}$ exists.
In what follows, $z_{0}$ is fixed  
thus  $|z|\sim |z+z_{0}|$ for large $|z|$.

We introduce some notation for  Sobolev spaces.

We denote the semi-classical $H^s$ norm by $\| w\|_{H^s_{sc}}^2=\int (1+h^2|\xi|^2)^s|\hat u(\xi)|^2d\xi$. 
Let $w$ be  a distribution on $\Omega$, we denote by $\| w\|_{\overline{ H}^s_{sc}(\Omega)}=\inf\{ \| \beta \|_{H^s_{sc}}, \text{ where }\beta_{|\Omega}=w  \}$. We recall that we denote  $D=-ih\d$,  and if $s$ is an integer the quantity $\sum_{|\alpha|\le s}\| D^\alpha w \|^2_{L^2(\Omega)}$ is equivalent to $\| w\|^2_{\overline{ H}^s_{sc}(\Omega)}$ uniformly with respect to $h$. 
When $h=1$ we denote the space by $\overline H^s(\Omega)$.

We apply the results of  \cite{Ro2013}. The estimates below are given by \cite[theorem~10]{Ro2013} for $k\ge 1$. 
The estimate on  $S_{12}$ and $S_{22} $ for $k=0$ are also given by \cite[theorem~10]{Ro2013}.
For $k=0$, \cite[Proposition~2.3]{Ro2013} gives the estimate on $S_{21}$  and  \cite[Lemma 2.1]{Ro2013} gives the estimate on $S_{11}$. The relation between $z$ and $h$ is $zh^2=\mu$, in particular $h^2|z|=1$.

\begin{align}
& S_{11} : \overline  H^{2k}_{sc}(\Omega) \to \overline  H^{2k+2}_{sc}(\Omega)  \text{ with } \| S_{11}  \|_{\overline  H^{2k}_{sc}(\Omega) \to \overline  H^{2k+2}_{sc}(\Omega)}\le C| z|^{-1}  \notag \\
  &S_{12} : \overline  H^{2k}_{sc}(\Omega) \to \overline  H^{2k+4}_{sc}(\Omega)   \text{ with } \|S_{12}  \|_{\overline  H^{2k}_{sc}(\Omega) \to \overline  H^{2k+4}_{sc}(\Omega)}\le C| z|^{-2}  \notag \\
& S_{21} : \overline  H^{2k}_{sc}(\Omega) \to \overline  H^{2k}_{sc}(\Omega)   \text{ with } \| S_{21}  \|_{\overline  H^{2k}_{sc}(\Omega) \to \overline  H^{2k}_{sc}(\Omega)}\le C  \notag\\
& S_{22} : \overline  H^{2k}_{sc}(\Omega) \to \overline  H^{2k+2}_{sc}(\Omega)   \text{ with } \|S_{22}  \|_{\overline  H^{2k}_{sc}(\Omega) \to \overline  H^{2k+2}_{sc}(\Omega)}\le C| z|^{-1}  \label{Estimations sur Sjk}
\end{align}

We denote by 
$$
\Lambda_z=
\begin{pmatrix}
\sqrt{| z|}&0\\ 0&1/\sqrt{| z|}.
\end{pmatrix}
$$
We remark that 
\begin{equation}\label{Form conjugaison par Lambda}
\Lambda_zA\Lambda_z^{-1}=
\begin{pmatrix}
A_{11} & |z|A_{12}\\
(1/|z|)A_{21}& A_{22}
\end{pmatrix} ,\text{ where }
A=
\begin{pmatrix}
A_{11} & A_{12}\\
A_{21}& A_{22}
\end{pmatrix}.
\end{equation}

We deduce from~\eqref{Estimations sur Sjk} that $ \Lambda_z S_{\omega_jz}\Lambda_z^{-1} : L^2(\Omega)\oplus L^2(\Omega) \to \overline  H^{2}_{sc}(\Omega)\oplus L^2(\Omega) $ with an operator  norm less than $C|z|^{-1}$ and $ \Lambda_z S_{\omega_jz}\Lambda_z^{-1} :\overline  H^{2k+2}_{sc}(\Omega)\oplus \overline  H^{2k}_{sc}(\Omega)\to H^{2k+4}_{sc}(\Omega)\oplus \overline  H^{2k+2}_{sc}(\Omega)$ with an  operator norm  less than  $C|z|^{-1}$.

As
\begin{equation}\label{Lien Tz et Sz}
\Lambda_z T_{z^p}\Lambda_z^{-1}=\prod_{j=1}^p  \Lambda_z S_{\omega_jz}\Lambda_z^{-1},
\end{equation}
we deduce that 
\begin{equation}\label{Operance de Tz conjugue}
\Lambda_z T_{z^p}\Lambda_z^{-1} :    L^2(\Omega)\oplus L^2(\Omega)\to \overline  H^{2p}_{sc}(\Omega)\oplus \overline  H^{2p-2}_{sc}(\Omega),
\end{equation}
with an  operator norm  less than  $C|z|^{-p}$. 

We can prove that $N(t)\le Ct^{n}$.  First we weaken~\eqref{Operance de Tz conjugue} to consider $\Lambda_z T_{z^p}\Lambda_z^{-1}$ as a map between $  L^2(\Omega)\oplus L^2(\Omega)\to \overline  H^{2p-2}_{sc}(\Omega)\oplus \overline  H^{2p-2}_{sc}(\Omega)$, with an  operator norm  less than  $C|z|^{-p}$.  
\newline  
As $\| v\|_{H^s}\le h^{-2s}\| v\|_{H^s_{sc}}$, we obtain
\begin{align}
&   \| \Lambda_z T_{z^p} \Lambda_z^{-1}\|_{L^2(\Omega)\oplus L^2(\Omega)\to \overline H^{2p-2}(\Omega)\oplus 
\overline  H^{2p-2}(\Omega)}\le C|z|^{-1}, \notag  \\
&\| \Lambda_z T_{z^p} \Lambda_z^{-1}\|_{L^2(\Omega)\oplus L^2(\Omega)\to L^2(\Omega)\oplus L^2(\Omega)}\le C|z|^{-p}. \label{Norm L2 de l'operateur conjugue}
\end{align} 
We can apply the theorem~13.5 in Agmon~\cite{Ag65}, that is, if $m>n/2$ we have
$$
\tno{T} \le C\|T\|^{n/(2m)}_m \|T\|^{1-n/(2m)}_0 ,
$$
where $\tno{T}$ is the Hilbert-Schmidt norm and $\|T\|_{m}$ is the operator norm of the map  $ L^2(\Omega)\oplus L^2(\Omega)\to \overline H^{m}(\Omega)\oplus \overline  H^{m}(\Omega)$. We apply this estimate with $m=2p-2>n/2$ and we have
$$
\tno{ \Lambda_z T_{z^p} \Lambda_z^{-1}}\lesssim |z|^{-\frac{n}{4(p-1)}} |z|^{-(1-\frac{n}{4(p-1)})p}= |z|^{-p+n/4}.
$$
We can follow the proof of Theorem~7 in~\cite{Ro2013}. If we denote by $\mu_{j}$ complex numbers such that $\mu_{j}^{-1}$ are eigenvalue of $S$ counted with multiplicity, then $\frac{1}{\mu_j^p-z^p} $ are the eigenvalues of $T_{z^p}$ and thus the eigenvalues of  $\Lambda_z T_{z^p} \Lambda_z^{-1}$. We obtain 
$$
\sum_j\frac{1}{|\mu_j^p-z^p|^2} \le\tno{ \Lambda_z T_{z^p} \Lambda_z^{-1}}^2\le C |z|^{-2p+n/2}.
$$
Let $\mu\in\C$ such that $|\mu|=1$, and $\omega_{j} \mu\not\in C_{e}\cup (0,\infty)$ for all $j=1,\cdots,p$. We take $z=t^{2}\mu$. If $|\mu_j|\le t^2$,  we have $ |\mu_j^p-z^p|\le 2t^{2p}$.
Then we have
$$
\sum_{|\mu_j|\le t^2}\frac 1{4t^{4p}}\le \sum_j\frac{1}{|\mu^p_j-z^p|^2} \le\tno{T_{z^p}}^2\le Ct^{-4p+n}.
$$
Then we obtain  $N(t)\le Ct^{n}$.

Now we prove Formula~\eqref{form entre serie et integrale}.
Estimate~\eqref{Operance de Tz conjugue} implies that $ \|  T_{11}\|_{L^2(\Omega)\to \overline  H^{2p}_{sc}(\Omega)}\le C|z|^{-p} $ and  as $\| v\|_{H^s}\le h^{-2s}\| v\|_{H^s_{sc}}$, we have
\begin{equation}\label{est H2p T11}
 \|  T_{11}\|_{L^2(\Omega)\to \overline  H^{2p}(\Omega)}\le C.
\end{equation}
To estimate the norm of $T_{22}$, we shall use that $S_{12}$ is a mapping from $L^2(\Omega)$ to $ \overline  H^{4}_{sc}(\Omega)$.
Actually if we take $g\in  L^2(\Omega)$, we have $ \Lambda_z S_{\omega_jz}\Lambda_z^{-1}(0,g)\in \overline  H^{4}_{sc}(\Omega)\oplus \overline  H^{2}_{sc}(\Omega)$ and $\| \Lambda_z S_{\omega_jz}\Lambda_z^{-1}(0,g)  \|_{\overline  H^{4}_{sc}(\Omega)\oplus \overline  H^{2}_{sc}(\Omega)}\le C  |z|^{-1}\| g \|_{L^2(\Omega)}$. We can repeat the previous argument for the $p-1$ other factors $  \Lambda_z S_{\omega_jz}\Lambda_z^{-1} $ and we obtain that 
$$
\|  \Lambda_z T_{z^p}\Lambda_z^{-1}(0,g) \|_{\overline  H^{2p+2}_{sc}(\Omega)\oplus \overline  H^{2p}_{sc}(\Omega)}\le C| z |^{-p}\| g \|_{L^2(\Omega)}.
$$
In particular this means that $\|T_{22}g\|_{\overline  H^{2p}_{sc}(\Omega)}\le C| z |^{-p}\| g \|_{L^2(\Omega)}$, and  
\begin{equation}\label{est H2p T22}
\|T_{22}g\|_{\overline  H^{2p}(\Omega)}\le C\| g \|_{L^2(\Omega)}. 
\end{equation}

By~\eqref{Norm L2 de l'operateur conjugue}, we have

\begin{equation}\label{est L2 operateurs diagonaux}
 \| T_{11} \|_{L^2(\Omega)\to L^2(\Omega)}+ \| T_{22}\|_{L^2(\Omega)\to L^2(\Omega)}\le  C|z|^{-p}.
\end{equation}

We can apply the theorem~13.9, Agmon~\cite{Ag65}. If $2p>n$ we have for $j=1,2$, $K_{jj}\in  \overline  H^{q}(\Omega\times\Omega )$ where $q=2p-[n/2]-1>n/2$. In particular the trace $K_{jj}(x,x)$ is well defined in $L^2(\Omega)$ and
\begin{align}
\left(
|  K_{jj}(x,x) |^2dx
\right)^{1/2}
&\le C\Big[
\left(
\| T_{jj} \|_{L^2(\Omega) \to \overline  H^{2p}(\Omega)}  +\| T_{jj}^* \| _{L^2(\Omega) \to \overline  H^{2p}(\Omega)} 
\right)^{n/(2p)}
\|T_{jj}  \|_{L^2(\Omega)\to L^2(\Omega)} ^{1-n/(2p)}  \notag \\
&\quad +\|T_{jj}  \|_{L^2(\Omega)\to L^2(\Omega)}  
\Big]   \label{Est. L2 du noyau par norme operateur}  
\end{align}

\begin{remark}\label{rem sur adjoint Sz}
The adjoint of $B_{z}$ is given by an analogous formula than~\eqref{eq : resolvent 1}. Indeed we find that the adjoint $B_z^*(p,q)=(g_1,g_2)$ is given by 
\begin{align*}
&- \Delta((1+\bar m)^{-1} p)-\bar z p=g_1\text{ in }\Omega   \\  
&-\Delta q-\bar zq+\bar m(1+\bar m)^{-1} p=g_2\text{ in }\Omega   \\
&q_{|\d\Omega}=\d_\nu q_{|\d\Omega}=0  \text{ on }\d\Omega.
\end{align*}
Using the relation between $R_{z}$ and $S_{z}$ (see~\eqref{Lien S et R}), we deduce that the adjoint of $S_{z}$ satisfies the same estimate than $S_{z}$ given in \eqref{Estimations sur Sjk}. By \eqref{Lien Tz et Sz}, the adjoint of $T_{z}$ satisfies the same properties than $T_{z}$ given in  \eqref{est H2p T11}, \eqref{est H2p T22} and \eqref{est L2 operateurs diagonaux}.
\end{remark}

 For $j=1,2$,\eqref{Est. L2 du noyau par norme operateur} implies, from \eqref{est H2p T11}, \eqref{est H2p T22} and \eqref{est L2 operateurs diagonaux}
\begin{equation}\label{Est. Noyau de Tzp}
\left(
\int_\Omega |K_{jj}(x,x)|^2dx
\right)^{1/2}\le C|z|^{-p+(n/2)}
\end{equation}
We recall that  $\lambda_j^{-1}$ are  the eigenvalues of $S$ counted with multiplicity 
the eigenvalues of $T$ are $\lambda_j^{-p}$. The indices are such that $|\lambda_j|\le |\lambda_{j+1}|$.
As $N(t)\le Ct^{n}$,
 this implies that $|\lambda_j|\ge C j^{2/n}$ where $C>0$. In particular $\sum 1/|\lambda_j|^p$ converges if $2p> n$.

By Theorem~12.17 in Agmon~\cite{Ag65}, there exists a constant $c\in\C$ such that
\begin{equation}\label{est trace}
\Tr(TT_{z^p})=\sum_{j\in\N} \frac{1}{(\mu_j^p-z^p)\mu_j^p}+c.
\end{equation}
We recall that the trace is defined in the theorem~12.20 in Agmon~\cite{Ag65} for an operator $Q=Q_1Q_2$ where $Q_1$ and $Q_2$ are Hilbert-Schmidt operators. Moreover if $K$ is the kernel of $Q$, $K(x,x)$ is definite for almost all $x$, we have $\int_{\Omega} |K(x,x) |dx< \infty $  and 
\begin{equation}\label{formule trace noyau}
\Tr(Q)=\int_{\Omega} K(x,x)dx.
\end{equation}

We remark as by assumption $\mu_j^p$ and $z^p$ are not in the same cone, we have $|\mu_j^p-z^p|\sim |\mu_j^p|+|z^p|$. Then 
$$
\frac{1}{|\mu_j^p-z^p||\mu_j^p|}\le \frac{C}{|\mu_j^p|} \text{ and } \frac{1}{|\mu_j^p-z^p||\mu_j^p|}\to 0 \text{ when } |z|\to +\infty.
$$
This implies that $\sum_{j\in\N} \frac{1}{(\mu_j^p-z^p)\mu_j^p}\to 0$ when $ |z|\to +\infty $.

In \cite{Ro2013} we  proved that $\tno{T_{z^p}}\le C | z|^{1-p+n/4}$ as $1-p+n/4<0$ thus $\tno{T_{z^p}}$ goes to 0 as $|z|$ goes to $+\infty$. We have $|\Tr(TT_{z^p})|\le \tno{T}\tno{T_{z^p}}$ goes to 0 as $|z|$ goes to $+\infty$. Then $c=0$ in~\eqref{est trace}. We obtain that 
$$
\Tr(z^pTT_{z^p})=\sum_{j\in\N}\left( \frac{1}{\mu_j^p-z^p}-\frac{1}{\mu_j^p}
 \right).
$$
We have $z^pTT_{z^p}=z^pT^2(I-z^pT)^{-1}=-T(I-z^pT)(I-z^pT)^{-1}+T(I-z^pT)^{-1}=T_{z^p}-T$. By Formula~\eqref{formule trace noyau} the trace of $z^pTT_{z^p}$ is given  by the integral of its kernel and, as the integral of the kernel of $T_{z^p}$ exists by \eqref{Est. Noyau de Tzp}, the integral of trace of kernel of $z^pTT_{z^p}-T_{z^p}$ does not depend of $z$,  we obtain
\begin{equation}\label{eq: noyau serie}
\int_\Omega K_{11}(x,x)dx+\int_\Omega K_{22}(x,x)dx=\sum_{j\in\N}\left( \frac{1}{\mu_j^p-z^p}
-\frac{1}{\mu_j^p}
 \right) +c
\end{equation}

By \eqref{Est. Noyau de Tzp}, for $j=1,2$, $\int_\Omega K_{jj}(x,x)dx$ goes to 0 when $|z|$ goes to 0 and as 
 $ \frac{1}{|\mu_j^p-z^p|}\le C \frac{1}{|\mu_j^p|}$ and $ \frac{1}{|\mu_j^p-z^p|}\to 0$ when $|z|\to+\infty$, we obtain $c=\sum_{j\in\N}
\frac{1}{\mu_j^p}$. Then \eqref{eq: noyau serie} gives the statement of  Lemma~\ref{Lem. Kernel trace}.

\end{subsection}

\begin{subsection}{Proof of Lemma \ref{Lem. trace symbol}}

We recall some facts on  pseudo-differential operators.
Let $a(x,\xi)$  be in  $\Con^\infty(\R^n\times\R^n)$ we say that $a$ is a symbol of order $m$ if for all $\alpha, \beta\in\N^n$, there exist $C_{\alpha, \beta}>0$,  such that
$$
|\d_x^\alpha\d_\xi^\beta a(x,\xi)|\le C_{\alpha, \beta}\est{\xi}^{m-|\beta|},
$$
where $\est{\xi}^2=1+|\xi|^2$. In particular a polynomial in $\xi$ of order $m$ with coefficients in $\Con^\infty(\R^n)$ with  bounded derivatives of all orders, is a symbol of order $m$.

With  a symbol we can associate an semi-classical operator by the following formula
$$
\Op(a)u=a(x,D)u= \frac{1}{(2\pi)^n}\int e^{ix\xi}a(x,h\xi)\hat u(\xi)d\xi= \frac{1}{(2h\pi)^n}\int e^{ix\xi/h}a(x,\xi)\hat u(\xi/h)d\xi.
$$
If $a(x,\xi)$ is a symbol of order $m$ (which can depend of $h$), if $a(x,\xi)=b(x,\xi)+hc(x,\xi)$ where $b(x,\xi)$ and $c(x,\xi)$ are symbols of order $m$, we call $b(x,\xi)$ the principal symbol of $a(x,\xi)$ which is definite modulo $h$.
This formula makes sense for $u\in\SSS(\R^n)$ and we can extend it to $u\in H^s$ for all $s$. For $a$, a symbol of order $m$, there exists $C>0$ such that for all $u\in H^s$,
$$
\| a(x,D)u\|_{H_{sc}^{s-m}}\le C\| u\|_{H_{sc}^{s}}. 
$$
In the following, when we use pseudo-differential operator we have always cut-off functions
supported in $\Omega$ in each side of the operator. We do not have to consider the action of pseudo-differential operator on $\overline{H}^{s}(\Omega)$ space as in~\cite{Ro2013}. 
 
 We begin with a description on $S_{z}$ in all compact set in $\Omega$.
 \begin{lemma}\label{lem description Sz}
  Let $\mu\in\C$ such that $|\mu|=1$ and we assume that 
    $\mu\not\in C_{e}\cup (0,+\infty)$
for $j=1,\cdots, p$. Let $z\in\C$ such that $z/\mu \in (0,\infty)$ is large enough.
   Let $ \theta$ and $\tilde\theta$ be functions in $\Con_0^\infty(\Omega)$ such that $\theta(x)=1$ if $x$ in the support of $\tilde \theta$, and $\theta(x)=1$ if $x$ in a compact subset of $\Omega$. Then we have 
 
\begin{equation}\label{Equation sur la parametrixe S}
\tilde\theta \Lambda_zS_z\Lambda_z^{-1}=|z|^{-1}\tilde\theta B \theta +|z|^{-1/2}W \theta\Lambda_zS_z\Lambda_z^{-1},
\end{equation}

where $W=\Lambda_zK\Lambda_z^{-1}$ and $W^{*} $ are  bounded on  $\overline{H}^s(\Omega)$ and the principal symbol of $B$ is given on the support of $\tilde\theta $ by 
\begin{equation*}
\begin{pmatrix}
 (a|\xi|^2-\mu)^{-1}&(a|\xi|^2-\mu)^{-1}V(|\xi|^2-\mu)^{-1}\\0&(|\xi|^2-\mu)^{-1}
\end{pmatrix} .
\end{equation*}
 \end{lemma}
\begin{prooff} 
 We apply the result proved in~\cite{Ro2013}. Let us recall the notations and the main results.

We multiply  Equations \eqref{eq : resolvent 1} by $h^2$, we denote by $\mu= h^2z$ where $\mu$ belongs to a bounded domain of $\C$, $a=1/(1+m)$ and $V=m/(1+m)$. We change $(f,g)$ in $(-f,-g)$.

We recall the assumption made on $m$, we have $m(x)\not=-1$ for all $x\in\overline \Omega$ and $m(x)\not=0$ for $x$ in a \nhd of $\d\Omega$.

Thus  following \eqref{eq :  resolvent 1}, we obtain the system
\begin{equation}\label{eq : resolvent sc}
\left\lbrace 
 	\begin{array}{ll}
 		&\big(-ah^2 \Delta -\mu\big) u-h^2Vv=h^2f\text{ in } \Omega\\
		&(-h^2\Delta -\mu)v=h^2g \text{ in } \Omega\\
		&u=\partial_\nu u=0 \text{ on } \partial \Omega.
	\end{array}\right. 
\end{equation}

In case where $n(x)=n_1(x)+n_2(x)/k$ we multiply Equations~\eqref{eq : resolvent 1 bis} by $h^2$, we denote by $\mu=-h^2k^2$, $a=1/(1+m_1)$, $V=m_1/(1+m_1) + h m_2/(\nu+\nu m_1)$ where $\nu=hk$. We change $(f,g)$ in $(-f,-g)$.
Thus following~\eqref{eq : resolvent 1 bis} we obtain the system
\begin{equation}\label{eq : resolvent sc bis}
\left\lbrace 
 	\begin{array}{ll}
 		&\big(-ah^2 \Delta+ h W_0 -\mu\big) u-h^2Vv=h^2f\text{ in } \Omega\\
		&(-h^2\Delta -\mu)v=h^2g \text{ in } \Omega\\
		&u=\partial_\nu u=0 \text{ on } \partial \Omega,
	\end{array}\right. 
\end{equation}
where $W_0=-\nu m_2/(1+m_2)$. In particular the principal semi-classical symbol of $-ah^2 \Delta+ h W_0 -\mu$ is $-a|\xi |^2-\mu$, the principal semi-classical symbol of $V$ is $m_1/(1+m_1)$. In what follows only the principal symbols of $-ah^2 \Delta+ h W_0 -\mu$ and 
$-h^2\Delta -\mu$
must be take account. For simplicity we write the proof for the system~\eqref{eq : resolvent sc}, the case of system~\eqref{eq : resolvent sc bis} may be treat following the same way.

Now we compute the symbol of the resolvent in $\Omega$.
 
 Let $\phi_0$, $\phi_1$ and $\phi_2\in\Con_0^\infty(\Omega)$ where $\phi_0\phi_1=\phi_0$ and $\phi_1\phi_2=\phi_1$. We take $\phi_{0}$  such that  $\phi_{0}=1$ on the support of $\tilde\theta$. Let $Q$ be a parametrix of $-ah^2\Delta-\mu$ such that $\phi_0 Q\phi_1(-ah^2\Delta-\mu )=\phi_0 -hK$ where $K$ is of order $-1$ and the principal symbol of $Q$ is  $(a|\xi|^2-\mu)^{-1}$. As
 $$
  \phi_1\big(-ah^2 \Delta -\mu\big)\phi_2 u=h^2V\phi_1 v+h^2 \phi_1 f\text{ in } \R^n,
 $$
 applying $\phi_{0} Q$ to this equation, we obtain
\begin{equation}\label{param. sur u fonction v}
 \phi_0 u=h^2   \phi_0
 Q \phi_1 f+h^2 \phi_0 Q\left( V\phi_1 v\right)+hK \phi_2 u.
\end{equation}

We apply the same method use above to compute  $v$  (see \cite[Section 2.2]{Ro2013}). Let $\phi_3\in\Con_0^\infty(\Omega)$ where $\phi_2\phi_3=\phi_2$  and $\theta =1$ on the support of $ \phi_{3}$. The choices of the $\phi_{j} $ are compatible with $\theta $ and $\tilde\theta$.
We have
$  \phi_1   \tilde Q\phi_2(-h^2\Delta-\mu)=\phi_1 -hK_{-1}$ where  $K_{-1} $ is of order $-1$ and the principal symbol of $\tilde Q$ is  $(|\xi|^2-\mu)^{-1}$.

We apply $\phi_{2}$ on the equation on $v$ in \eqref{eq : resolvent sc}, we have
$$
\phi_{2}(-h^2\Delta -\mu)\phi_{3}v=h^2\phi_{2}g.
$$
Applying the parametrix $\phi_{1}\tilde Q$, we have
\begin{equation}\label{Param sur v}
\phi_1v=h^2 \phi_1\tilde Q \phi_2 g+hK_{-1}\phi_3 v.
\end{equation}
With this equation and \eqref{param. sur u fonction v} we obtain

\begin{equation}\label{Param sur u}
\phi_0 u=h^2\phi_0Q \phi_1 f+h^4 \phi_{0}Q\left( V\phi_1 \tilde Q\phi_2 g\right)+hK \phi_2 u +h^3\tilde K_{-1}\phi_3 v,
\end{equation}
where $\tilde K_{-1}$ is an operator of order $-1$ .

Let 
\begin{equation*}
A=
\begin{pmatrix}
h^2 Q\phi_{1}&h^4QV\phi_{1}\tilde Q \phi_{2}\\0&h^2\phi_{1}\tilde Q\phi_{2}
\end{pmatrix} 
\text{ and }
K=
\begin{pmatrix}
 K_{11} &h^2K_{12}\\ 0& K_{22}
\end{pmatrix},
\end{equation*}
where $K_{11}=  K\phi_{2}$, $K_{12}= \tilde K_{-1} $ and $K_{22} = K_{-1}\phi_{2}  $.  In particular   $K_{jk}$ is bounded on $\overline{H}^s_{sc}(\Omega)$ for all $s\ge0$. Indeed, all the operators contain cut-off thus $K_{jk}\phi_{3}u$ is compactly supported in $\Omega$ if $u\in  \overline{H}^s_{sc}(\Omega)$.

We recall that $S=R_{z_{0}} $ and by~\eqref{Lien S et R}, $S_z=R_{z_{0}+z}$,

if the resolvent $R_{z_{0}+z}$ exists.
In what follows, $z_{0}$ is fixed and we have the relation  $z=-\mu/h^{2}-z_{0}$. In particular $|z|^{-1/2}\sim h$ for large $|z|$. 
With these relations we have,  following \eqref{Param sur v} and \eqref{Param sur u},
 $S_z(f,g)=(u,v)$, and
 
\begin{equation*}
\phi_0 S_z=\phi_0 A\phi_3 +|z|^{-1/2}K\phi_3 S_z.
\end{equation*} 

Thus we can write
\begin{equation}\label{eq. Equation sur la parametrixe S}
\phi_0 \Lambda_zS_z\Lambda_z^{-1}=|z|^{-1}\phi_0 B\phi_3 +|z|^{-1/2}W\phi_3 \Lambda_zS_z\Lambda_z^{-1},
\end{equation}

where $W=\Lambda_zK\Lambda_z^{-1}$ is bounded on  $\overline{H}^s(\Omega)$ and the principal symbol of $B$ is given on the support of $\phi_{0}$ by 
\begin{equation*}
\begin{pmatrix}
 (a|\xi|^2-\mu)^{-1}&(a|\xi|^2-\mu)^{-1}V(|\xi|^2-\mu)^{-1}\\0&(|\xi|^2-\mu)^{-1}
\end{pmatrix} .
\end{equation*}
As $\tilde\theta\phi_{0}=\tilde\theta$ and $\phi_{3}\theta=\phi_{3}$, \eqref{eq. Equation sur la parametrixe S} gives~\eqref{Equation sur la parametrixe S}. As $K$ is a semi-classical pseudo-differential operator, $W^{*} $ is also bounded on $\overline{H}^s(\Omega)$
\end{prooff}
 
\begin{remark}
Formula~\eqref{Equation sur la parametrixe S} does not give a description on the operator $S_{z}$
 in $\Omega$. It gives only $S_{z} $ in all compact in $\Omega $. In the proof below we need also 
 estimates on $S_{z}$ up the boundary given in~\eqref{Estimations sur Sjk} to absorb the error terms.
\end{remark}

Lemma~\ref{lem description Sz} gives the principal symbol of $S_{z}$, the following lemma gives the principal symbol of 
$T_{z^{p}}$.  
\begin{lemma}\label{lem: formule de Tzp}
Let $p\in \N\setminus\{ 0 \}$. Let $\varphi_{0}$ and $\varphi_{1}$ be functions in $\Con_0^\infty(\Omega)$ such that $\varphi_{0}=1$ on a \nhd of $\supp \varphi_{1}$.
\begin{equation}\label{Formule pour Tzp}
\varphi_1 \Lambda_z T_{z^p} \Lambda_z^{-1}\varphi_0=|z|^{-p}\varphi_1B_p\varphi_0+|z|^{-p-1/2}\varphi_1R_p\varphi_0,
\end{equation}
where  $\varphi_1R_p\varphi_0$ satisfies  the following property if we denote by 
$
\varphi_1R_p \varphi_0 =
\begin{pmatrix}
R_p^{11}&R_p^{12}\\R_p^{21}&R_p^{22}
\end{pmatrix}
$

\begin{align}
&R_p^{11} : L^2(\Omega) \to  \overline H^{2p}_{sc}(\Omega)  \notag \\
&R_p^{12} : L^2(\Omega) \to   \overline H^{2p+2}_{sc}(\Omega)  \notag \\
&R_p^{21} :  L^2(\Omega) \to  \overline H^{2p-2}_{sc}(\Omega)   \notag\\
&R_p^{22} :  L^2(\Omega) \to  \overline H^{2p}_{sc}(\Omega) , \label{Estimations sur Rjk}
\end{align}
where the norm operator are uniformly  bounded with respect $h$, 

 and the principal symbol of $B_{p}$ is 
\begin{equation}\label{form Bp lemme}
\begin{pmatrix}
(a^p|\xi|^{2p}-\mu^p)^{-1} & Q_{-2p-2}(x,\xi)\\ 0 &(|\xi|^{2p}-\mu^p)^{-1}
\end{pmatrix}
\end{equation}
Moreover the adjoint of $R^{jq}_{p}$ satisfies~\eqref{Estimations sur Rjk} where the norm operator are uniformly  bounded with respect $h$.
\end{lemma}
\begin{prooff}
We argue by induction on $k$ and for that we must introduce a sequence of cut-off functions. Let $\chi_k$ and $\tilde\chi_k$ be cut-off functions such that $\tilde\chi_k \tilde\chi_{k+1} = \tilde\chi_{k+1} $, $\chi_k \tilde\chi_{k} = \tilde\chi_{k} $,   $\tilde\chi_k \chi_{k+1} = \chi_{k+1} $. We can assume that $\tilde\chi_{0}=1$ on the support of 
$\varphi_{0}$ and and $\tilde\chi_{p}=1$ on the support of $\varphi_{1}$. We can apply Formula~\eqref{Equation sur la parametrixe S} where $\tilde\theta$ is replaced by $\tilde\chi_k$ and $\theta$ by $\chi_k$. We have 
\begin{equation}\label{Equation sur la parametrixe S dependant de k}
\tilde\chi_k \Lambda_zS_z\Lambda_z^{-1}=|z|^{-1}\tilde\chi_k B\chi_k + |z|^{-1/2}W\chi_k \Lambda_zS_z\Lambda_z^{-1},
\end{equation}
where $W$ and $W^{*}$ are bounded on $\overline{H}^s(\Omega)$.

We prove by recurrence  the following formula
\begin{equation}\label{Eq. rec. sur pd pour Tz}
\tilde\chi_k\prod_{j=1}^k  \Lambda_z S_{\omega_jz}\Lambda_z^{-1}\tilde\chi_0=|z|^{-k}\tilde\chi_kB_k\tilde\chi_0+|z|^{-k-1/2} 
R_k\tilde\chi_0,
\end{equation}
where the semi-classical principal symbol of $B_k$ is given by
\begin{equation}\label{Formule pour Bk}
B_k=
\begin{pmatrix}
 \prod _{j=1}^k(a|\xi|^2-\omega_j\mu)^{-1}&  Q_{-2k-2}(x,\xi) \\ 0&  \prod _{j=1}^k(|\xi|^2-\omega_j\mu)^{-1}
\end{pmatrix},
\end{equation}
where $Q_{-2k-2}(x,\xi)$ is a symbol of order $-2k-2$.
The operators $R_k ^{jq}$  and their adjoints satisfy Estimates~\eqref{Estimations sur Rjk} with $p=k$.

For $k=1$, Formula~\eqref{Formule pour Bk} and 
properties~\eqref{Estimations sur Rjk} for $R_1 ^{jq} $ 	and their adjoints follow from~\eqref{Form conjugaison par Lambda}, \eqref{Estimations sur Sjk}, Remark~\ref{rem sur adjoint Sz} and \eqref{Equation sur la parametrixe S dependant de k}.

If Formula~\eqref{Eq. rec. sur pd pour Tz} is true for $k$ we have 
\begin{equation*}
\tilde\chi_{k+1}\prod_{j=1}^{k +1} \Lambda_z S_{\omega_jz}\Lambda_z^{-1}\tilde\chi_0= L_1+L_2,
\end{equation*}
where
\begin{align*}
L_1 & =    \tilde\chi_{k+1} \Lambda_z S_{\omega_{k+1}z}\Lambda_z^{-1} \tilde\chi_k\prod_{j=1}^k  \Lambda_z S_{\omega_jz}\Lambda_z^{-1}\tilde\chi_0\\
L_2  & =  \tilde\chi_{k+1} \Lambda_z S_{\omega_{k+1}z}\Lambda_z^{-1}(1- \tilde\chi_k)\prod_{j=1}^k  \Lambda_z S_{\omega_jz}\Lambda_z^{-1}\tilde\chi_0.
\end{align*}
By~\eqref{Equation sur la parametrixe S dependant de k} for $k+1$ and \eqref{Eq. rec. sur pd pour Tz} for $k$, we have
\begin{equation*}
L_1= \left( |z|^{-1} \tilde\chi_{k+1}B\chi_{k+1} +|z|^{-1/2}W \chi_{k+1} \Lambda_zS_z\Lambda_z^{-1}\right)
\left( |z|^{-k}\tilde\chi_kB_k\tilde\chi_0+|z|^{-k-1/2}  
R_k \tilde\chi_0\right).
\end{equation*}
The term $ |z|^{-1-k} \tilde\chi_{k+1}B\chi_{k+1} \tilde\chi_kB_k\tilde\chi_0$ gives the first right hand side term of~\eqref{Eq. rec. sur pd pour Tz}, where $B_{k+1}=B\chi_{k+1}\tilde\chi_{k}B_{k}$ and the principal symbol is given by  Formula~\eqref{Formule pour Bk} on the support of $\tilde\chi_{k+1}$. 

The three other terms have the form of $R_{k+1}$ and satisfy  Estimates~\eqref{Estimations sur Rjk}. Indeed, the power of $|z|$ is obtained as the operator norm of $\Lambda_zS_z\Lambda_z^{-1}$ is bounded by $|z|^{-1}$. To prove the mapping between the $H^{s}_{sc}$, we denote by $A_{q}$ a generic operator of order $q$ mapping 
$H^{s}$ to $H^{s-q}$. We check that 
\begin{equation}\label{composition operateur inhomogene}
\begin{pmatrix}
A_{-2} & A_{-4}\\ A_{0} &A_{-2}
\end{pmatrix}
.
\begin{pmatrix}
A_{-2k} & A_{-2k-2}\\ A_{-2k+2} &A_{-2k}
\end{pmatrix}=
\begin{pmatrix}
A_{-2k-2} & A_{-2k-4}\\ A_{-2k} &A_{-2k-2}
\end{pmatrix}.
\end{equation}
The properties on adjoints follow from the recurrence assumptions on $R_{k}$, the properties on $W^{*}$ and Remark~\ref{rem sur adjoint Sz}.

By \eqref{Equation sur la parametrixe S dependant de k}  for $k+1$ and as $\chi_{k+1}(1-\tilde\chi_k)=0$, we have

\begin{equation*}
L_2= |z|^{-1/2}W\chi_k \Lambda_zS_z\Lambda_z^{-1}     (1- \tilde\chi_k)\prod_{j=1}^k  \Lambda_z S_{\omega_jz}\Lambda_z^{-1}\tilde\chi_0.
\end{equation*}
By \eqref{Form conjugaison par Lambda} and 
\eqref{Estimations sur Sjk} the operator norm of this term is $|z|^{-k-1/2}$. The proof that $L_{2}$ satisfies Estimates~\eqref{Estimations sur Rjk} for $p=k+1$ is obtained by~\eqref{composition operateur inhomogene}. The properties on the adjoint of $L_{2}$ follow from Remark~\ref{rem sur adjoint Sz} and the properties on $W^{*}$.

From~\eqref{Formule pour Bk} for $k=p$, and as  $z^p-\mu^p=\prod_{j=1}^p(z-\omega_j\mu)$, we obtain that $ \prod _{j=1}^p(a|\xi|^2-\omega_j\mu)^{-1}=(a^{p}|\xi|^{2p}-\mu^{p})^{-1}  $. This gives the diagonal terms of the symbol of $B_{p}$ in Formula~\eqref{form Bp lemme}. 

\end{prooff}

\medskip

Now we can finish the proof of Lemma~\ref{Lem. trace symbol}.
We take $\varphi_0$ such that $\varphi_0(x)=1$ is $d(x, \R^n\setminus \Omega)\ge 2\delta$ and  $\varphi_0(x)=0$ is $d(x, \R^n\setminus \Omega)\le \delta$. We take $\varphi_1$ such that $\varphi_1(x)=1$ is $d(x, \R^n\setminus \Omega)\ge 4\delta$ and  $\varphi_1(x)=0$ is $d(x, \R^n\setminus \Omega)\le 3\delta$.

With the notation of \eqref{Formule pour Tzp} we deduce from \eqref{Estimations sur Rjk} that for $j=1,2$ we have 
$$
\| |z|^{-p-1/2} \varphi_1R_p^{jj}\varphi_0\|_{L^2(\Omega)\to \overline  H^{2p}_{sc}(\Omega)}\lesssim |z|^{-p-1/2}.
$$
This implies as for $T_{jj}$ in \eqref{est H2p T22}  and \eqref{est L2 operateurs diagonaux} that 
\begin{equation*}
\| |z|^{-p-1/2} \varphi_1R_p^{jj}\varphi_0\|_{L^2(\Omega)\to L^2(\Omega)}\lesssim |z|^{-p-1/2} \text{ and }
\| |z|^{-p-1/2} \varphi_1R_p^{jj}\varphi_0\|_{L^2(\Omega)\to \overline  H^{2p}(\Omega)}\lesssim |z|^{-1/2}.
\end{equation*}
By Formula~\eqref{Est. L2 du noyau par norme operateur} 	applied to the kernel of $|z|^{-p-1/2} \varphi_1R_p^{jj}\varphi_0$ denoted by $K_{jj}^R$ and the properties on $R_{p}$ and its adjoint given in Lemma~\ref{lem description Sz}, we obtain
\begin{equation}\label{Est. norm L2 reste Rp}
\left(
\int_\Omega |K_{jj}^R(x,x)|^2dx
\right)^{1/2}\le C|z|^{-p+(n/2)-1/2}.
\end{equation}

By the principal symbol of $B_{p}$ given in Lemma~\ref{lem: formule de Tzp} we can compute the integral of the kernel $K_{jj}(x,x)$.
Denoting by $b(x,\xi)$ either $(a^p|\xi|^{2p}-\mu^p)^{-1} $ or $(|\xi|^{2p}-\mu^p)^{-1} $, the kernel of a diagonal term is given by $(2\pi h)^{-n}\int e^{i(x-y)\xi/h}\varphi_1(x)b(x,\xi)\varphi_0(y)d\xi$ and this integral make sense if $p> n/2$. Denoting by $ K_{jj}^{B_{p}}(x,y)$ the diagonal terms of the kernel of $|z|^{-p} \varphi_{1}B_{p}\varphi_{0}$, we obtain
\begin{multline}\label{Debut lien trace symbole}
\int_\Omega K_{11}^{B_{p}}(x,x)dx+\int_\Omega K_{22}^{B_{p}}(x,x)dx
\\=(2\pi)^{-n}|z|^{-p+n/2}\int\int\varphi_1(x)\left((a^p|\xi|^{2p}-\mu^p)^{-1} +(|\xi|^{2p}-\mu^p)^{-1}
\right)d\xi dx +O(|z|^{-p+(n+1)/2}),
\end{multline}
where the error term $O(|z|^{(n+1)/2})$ is given by the lower order terms in the symbolic calculus.
From \eqref{Formule pour Tzp}, \eqref{Est. norm L2 reste Rp} and \eqref{Debut lien trace symbole}, we deduce 
\begin{multline*}
\int \varphi_{0}(x)(K_{11}(x,x)+K_{22}(x,x))dx
\\
=(2\pi)^{-n}|z|^{-p+n/2}\int\int\varphi_1(x)\left((a^p|\xi|^{2p}-\mu^p)^{-1} +(|\xi|^{2p}-\mu^p)^{-1}
\right)d\xi dx +O(|z|^{-p+(n+1)/2}).
\end{multline*}

Now we can write
\begin{equation*}  
\Lambda_z T_{z^p}\Lambda_z ^{-1}=\varphi_1  \Lambda_z T_{z^p}\Lambda_z ^{-1}   \varphi_0+\varphi_1  \Lambda_z T_{z^p}\Lambda_z ^{-1}  (1- \varphi_0)+  (1-\varphi_1 ) \Lambda_z T_{z^p}\Lambda_z ^{-1}.
\end{equation*}
If $K(x,y)$ is the kernel of $\Lambda_z T_{z^p}\Lambda_z ^{-1}$, the kernel of the left hand side terms are respectively, 
$$
\varphi_{1}(x)K(x,y)  \varphi_0(y),  \ \varphi_{1}(x)K(x,y)(1-  \varphi_0(y)) \text{ and  }(1-\varphi_{1}(x))K(x,y) .
$$ 
In particular, from the properties of the supports of $\varphi_{j}$, we have  $\varphi_{1}(x)K(x,x)(1-  \varphi_0(x))=0$.
Let $F(x)=K_{11}(x,x)+K_{22}(x,x)$, to show that 
\begin{equation}\label{maj trace avec delta}
\left|\int_{\Omega}   |z|^{p-n/2} F(x)dx
-(2\pi)^{-n}\int_\Omega\int\left((a^p|\xi|^{2p}-\mu^p)^{-1} +(|\xi|^{2p}-\mu^p)^{-1}
\right) d\xi dx  \right|\le C\delta+C_\delta | z |^{-1/2},
\end{equation}
We shall prove
\begin{multline}\label{Est. Trace moins int. symb.}
\left| |z|^{p-n/2}  \varphi_1(x) F(x)  dx
-(2\pi)^{-n}\int_\Omega\int\left((a^p|\xi|^{2p}-\mu^p)^{-1} +(|\xi|^{2p}-\mu^p)^{-1}
\right) d\xi dx  \right|
\\
\le C\delta+C_\delta | z |^{-1/2} ,
\end{multline}
and
\begin{equation}\label{Est. rest. pour Trace}
    |z|^{p-n/2}   \left|
  \int_\Omega(1-  \varphi_1(x))F(x)dx
    \right|\le C \delta^{1/2}.
\end{equation}
Obviously \eqref{Est. Trace moins int. symb.} and \eqref{Est. rest. pour Trace}  imply  \eqref{maj trace avec delta}, and \eqref{maj trace avec delta} implies Lemma~\ref{Lem. trace symbol}.
To prove \eqref{Est. Trace moins int. symb.}, we apply \eqref{Debut lien trace symbole}. We have 
\begin{multline*}
\int_{\{ x,\ d(x,\R^n\setminus\Omega)\le 2\delta \}}\int\left|(a^p|\xi|^{2p}-\mu^p)^{-1} +(|\xi|^{2p}-\mu^p)^{-1}
\right|d\xi dx
\\
\le C \int_{\{ x,\ d(x,\R^n\setminus\Omega)\le 2\delta \}}\int \est\xi^{-2p}d\xi dx\le C\delta,
\end{multline*}
indeed $|\xi|^{2p}$ and  $a^{p}|\xi|^{2p}$  are not in the same cone as $\mu^{p} $  by assumption. 
Thus the principal term in  \eqref{Est. Trace moins int. symb.}  given by the principal term from \eqref{Debut lien trace symbole}, can be estimate by $C\delta$ and the error terms from \eqref{Debut lien trace symbole} can be estimate by $C_{\delta}|z|^{-1/2}$.

To prove  \eqref{Est. rest. pour Trace}, using   \eqref{Est. Noyau de Tzp} we obtain
\begin{equation*}
   \left|
  \int_\Omega(1-  \varphi_1(x))K_{jj}(x)dx    \right|\le 
   \left(
\int_\Omega
(K_{jj}(x,x))^2 dx
\right)^{1/2}
\left(
\int_{\{x, d(x,\R^n\setminus\Omega)\le 4\delta\}} dx
\right)^{1/2}\le C|z|^{-p+n/2}\delta^{1/2},
\end{equation*}
which is Estimate~\eqref{Est. rest. pour Trace}.

\end{subsection}

\end{section}

\begin{section}{Proof of Theorem \ref{th. count. funct.}}

We shall apply the following tauberian Theorem cited in Agmon~\cite[th. 14.5]{Ag65}, The proof is given in Karamata~\cite{Ka31}.

\begin{theorem}[Tauberian Theorem]\label{Th. Tauberian}
Let $\sigma(\lambda)$ be a non decreasing function for $\lambda > 0$,
let $ 0 < a < 1$, let  $ \alpha$ be a non-negative number, and suppose that as $t\to +\infty$,
$$
\int_0^{+\infty}\frac{d\sigma(\lambda)}{\lambda+t}=\alpha t^{a-1}+o(t^{a-1}).
$$

Then as $\lambda\to +\infty$
$$
\sigma(\lambda)=\alpha \frac{ \sin \pi a}{ \pi a}\lambda^a+o(\lambda^a).
$$

\end{theorem}

Theorem~\ref{th. count. funct.} is implied by  tauberian theorem and Theorem~\ref{th. trace formula} which gives
\begin{align}
&\sum_{j\in\N} \frac{1}{\mu_j^p-z^p}=A|z|^{-p+n/2} +o(|z|^{-p+n/2}),  \notag\\
&\text{where } A=(2\pi)^{-n}\int_\Omega\int\left((a^p|\xi|^{2p}-\mu^p)^{-1} +(|\xi|^{2p}-\mu^p)^{-1}
\right) d\xi dx, \text{ with } \mu=z/|z|.  \label{definition de A}
\end{align}
We take $z=t^{1/p}e^{i\pi/p}$, where $t>0$, we obtain

\begin{equation}\label{Asymptotic trace kernel}
\sum_{j\in\N} \frac{1}{\mu_j^p+t}=At^{-1+n/2p} +o(t^{-1+n/2p}).
\end{equation}
Let $\mu_j=\delta_j+i\nu_j$ where $\delta_j$ and $\nu_j$ are real. Let $A(t)=\sum_{j\in\N} \frac{1}{\delta_j^p+t}$.

By assumption (see Remark~\ref{rem vp pres du reel}),  for all $\eps>0$, there exists $R_{\eps}>0$, such that if $\delta_j\ge R_{\eps} $ then $|\nu_j|\le \eps \delta_j$.

We have for $j$ large enough such that $|\nu_j|\le \eps \delta_j$,
\begin{equation}\label{Est. diff. real and imaginary}
\left| \frac1{\delta_j^p+t} -\frac{1}{\mu_j^p+t} \right|\le C\eps\frac1{\delta_j^p+t},
\end{equation}
where $C$ depends only on $p$.  

Denoting $N_\eps$ such that for all $j\ge N_\eps$, then $|\nu_j|\le \eps \delta_j$.
We have 
\begin{equation}\label{Sum first terms}
\sum_{j=1}^{N_\eps}\frac{1}{|\mu_j^p+t|}  \le C \sum_{j=1}^{N_\eps}\frac{1}{\delta_j^p+t}\le CN_\eps/t  .
\end{equation}

We deduce from \eqref{Asymptotic trace kernel},  \eqref{Est. diff. real and imaginary} and \eqref{Sum first terms} that 
\begin{equation}\label{form pour A(t)}
A(t)+O(N_\eps/t)+O(\eps A(t))=At^{-1+n/2p} +o(t^{-1+n/2p}).
\end{equation}
For $\eps$ small enough we deduce there exists $C>0$ such that 
$$
(1/C)t^{-1+n/2p}\le A(t)\le Ct^{-1+n/2p}.
$$
Using this in~\eqref{form pour A(t)} we obtain
\begin{equation*}
\sum_{j\in\N} \frac{1}{\delta_j^p+t}=At^{-1+n/2p} +o(t^{-1+n/2p}).
\end{equation*}

Denoting $\sigma(\lambda)=\sharp \{j\in\N,\  \delta_j^p\le\lambda \}$ where the number is counted with multiplicity. We have 
\begin{equation*}
\sum_{j\in\N} \frac{1}{\delta_j^p+t}=\int_0^{+\infty}\frac{d\sigma(\lambda)}{\lambda+t}=At^{-1+n/2p} +o(t^{-1+n/2p}).
\end{equation*}
By tauberian Theorem 
we obtain 
\begin{equation*}
\sharp \{j\in\N,\  \delta_j^p\le\lambda \}=2p A \frac{ \sin \big(  n\pi/(2p)\big) }{n \pi }\lambda^{n/(2p)}+o(\lambda^{n/(2p)}),
\end{equation*}
which is equivalent to 
\begin{equation}\label{form pour N(t)}
\sharp \{j\in\N,\  \delta_j\le t^2 \}=2p A \frac{ \sin \big(  n\pi/(2p)\big) }{n \pi }t^{n}+o(t^{n}).
\end{equation}

We recall that $\mu=z/|z|$ and $z=t^{1/p}e^{-i\pi/p}$ thus $\mu=-1$. To compute $A$ from  \eqref{definition de A}, we must compute integral as $\int (b^p|\xi|^{2p}+1)^{-1}d\xi$, where $b>0$. We have 
\begin{align*}
\int (b^p|\xi|^{2p}+1)^{-1}d\xi  & =  n\text{Vol}(B_1)\int_0^{+\infty}r^{n-1}(b^pr^{2p}+1)^{-1}dr  \\
     &  = n(2p)^{-1 } b^{-n/2}\text{Vol}(B_1)\int_0^{+\infty}\sigma^{n/(2p)-1}(\sigma+1)^{-1}d\sigma  \\
     &  =   n(2p)^{-1 } b^{-n/2}\text{Vol}(B_1) \pi \sin^{-1}(n\pi/(2p))  ,
\end{align*}
where the last integral is computed by residue theorem (see Cartan~\cite[p 107]{Ca}).

 This gives that 
\begin{align*}
2p A \frac{ \sin \big(  n\pi/(2p)\big) }{n \pi }=(2\pi)^{-n}\text{Vol}(B_1) \int_\Omega(a^{-n/2}(x) +1)dx=\alpha.
\end{align*}
Now to prove the statement of Theorem~\ref{th. count. funct.}, we must prove that  $\alpha= \lim t^{-n}\sharp \{j\in\N,\  \delta_j\le t^{2} \}=\lim  t^{-n}\sharp \{j\in\N,\  |\mu_{j}| \le t^{2} \}$ when $t$ goes to $\infty$. Except for a finite number of values, $\delta_{j}>0$ and as $\delta_{j}\le |\mu_{j}|$ there exists $C>0$ such that
$$
 \sharp \{j\in\N,\  |\mu_{j}| \le t^{2} \}  \le    \sharp \{j\in\N,\  \delta_j\le t^{2} \}+ C.
$$
 Thus 
\begin{equation}\label{Majoration de fct comptage}
\limsup_{t\to \infty}  t^{-n}\sharp \{j\in\N,\   |\mu_{j}|\le t^{2} \}\le \alpha.
\end{equation}
For all $\eps>0$, there exists $J_{\eps}$ such that for all $j\ge J_{\eps}$, $|\nu_{j}|\le \eps \delta_{j}$ then $|\mu_{j}|\le (1+\eps)\delta_{j}$. Thus there exist $C_{\eps}>0$ such that 
$$
  \sharp \{j\in\N,\  \delta_j\le t^{2} \}\le   \sharp \{j\in\N,\  (1+\eps)^{-1}|\mu_{j}| \le t^{2} \}   +C_{\eps},
$$
which is equivalent to 
$$
  \sharp \{j\in\N,\  \delta_j\le (1+\eps)^{-1}  t^{2} \}\le   \sharp \{j\in\N,\  |\mu_{j}| \le t^{2} \}   +C_{\eps}.
$$
We obtain
$$
 (1+\eps)^{-n/2} \alpha\le \liminf_{t\to\infty}t^{-n} \sharp \{j\in\N,\  |\mu_{j}| \le t^{2} \} .
$$
As this estimate is true for all $\eps>0$ and from~\eqref{Majoration de fct comptage}  we have 
$$
\alpha\le  \liminf_{t\to\infty}t^{-n} \sharp \{j\in\N,\  |\mu_{j}| \le t^{2} \}\le \limsup_{t\to \infty}  t^{-n}\sharp \{j\in\N,\   |\mu_{j}|\le t^{2} \}\le \alpha.
$$
This is the statement of Theorem~\ref{th. count. funct.}.

\end{section}

\end{document}